\documentclass{article}
\usepackage{amsmath}
\usepackage{amsfonts}

\newtheorem{theorem}{Theorem}

\newtheorem{example}[theorem]{Example}

\input{tcilatex}
\begin{document}

\title{A New Algorithm for the Inverse of Periodic $k$ Banded and Periodic
Anti $k$ banded Matrices}
\author{M. Ya\c{s}ar\thanks{%
e-mail: meral.yasar@hotmail.com} \& D. Bozkurt\thanks{%
e-mail: dbozkurt@selcuk.edu.tr} \and Department of Mathematics, Nigde
University and \and Department of Mathematics, Selcuk University}
\maketitle

\begin{abstract}
In this study, an algorithm for computing the inverse of periodic $k$ banded
matrices , which are needed for solving the differential\linebreak equations
by using the finite differences, the solution of partial differential
equations and the solution of boundary value problems is obtained and the
inverses of periodic anti $k$ banded matrices are computed. In addition,
the\linebreak determinant of these type of matrices and the solution of
linear systems having these coefficient matrices are investigated. When
obtaining this\linebreak algorithm, the $LU$ factorization is used.The
algorithm is implementable to the CAS (Computer Algebra Systems) such as
Maple and Mathematica.
\end{abstract}


\section{Introduction}

The $n\times n$ periodic $k$ banded and anti $k$ banded matrices take the
following form respectively:%
\begin{eqnarray}
M &=&\left[ 
\begin{array}{ccccc}
a_{11} & a_{12} & \cdots & a_{1,\frac{k+1}{2}} & 0 \\ 
a_{21} & a_{22} & \cdots & a_{2,\frac{k+1}{2}} & a_{2,\frac{k+3}{2}} \\ 
a_{31} & a_{32} & \cdots & a_{3,\frac{k+1}{2}} & a_{3,\frac{k+3}{2}} \\ 
\vdots & \vdots & \ddots & \vdots & \vdots \\ 
a_{\frac{k+1}{2},1} & a_{\frac{k+1}{2},2} & \cdots & a_{\frac{k+1}{2},\frac{%
k+1}{2}} & a_{\frac{k+1}{2},\frac{k+3}{2}} \\ 
0 & a_{\frac{k+3}{2},2} & \cdots & a_{\frac{k+3}{2},\frac{k+1}{2}} & a_{%
\frac{k+3}{2},\frac{k+3}{2}} \\ 
\vdots & \vdots &  & \vdots & \vdots \\ 
0 & 0 & \cdots & a_{n-\frac{k-1}{2},\frac{k+1}{2}} & a_{n-\frac{k-1}{2},%
\frac{k+3}{2}} \\ 
\vdots & \vdots &  & \vdots & \vdots \\ 
0 & 0 & \cdots & 0 & 0 \\ 
a_{n1} & 0 & \cdots & 0 & 0%
\end{array}%
\right.  \notag \\
&&\ \ \ \ \ \ \ \ \ \ \ \ \ \ \ \left. 
\begin{array}{cccc}
0 & \cdots & 0 & a_{1n} \\ 
0 & \cdots & 0 & 0 \\ 
a_{3,\frac{k+5}{2}} & \cdots & 0 & 0 \\ 
\vdots &  & \vdots & \vdots \\ 
a_{\frac{k+1}{2},\frac{k+5}{2}} & \cdots & 0 & 0 \\ 
a_{\frac{k+3}{2},\frac{k+5}{2}} & \cdots & 0 & 0 \\ 
\ddots &  & \vdots & \vdots \\ 
a_{n-\frac{k-1}{2},\frac{k+5}{2}} & \cdots & a_{n-\frac{k-1}{2},n-1} & a_{n-%
\frac{k-1}{2},n} \\ 
\vdots & \ddots & \vdots & \vdots \\ 
a_{n-1,n-\frac{k+1}{2}} & \cdots & a_{n-1,n-1} & a_{n-1,n} \\ 
0 & \cdots & a_{n,n-1} & a_{n,n}%
\end{array}%
\right]
\end{eqnarray}%
\begin{eqnarray}
N &=&\left[ 
\begin{array}{cccccc}
a_{11} & 0 & \cdots & 0 & 0 & 0 \\ 
0 & 0 & \cdots & 0 & 0 & a_{2,n-\frac{k+1}{2}} \\ 
0 & 0 & \cdots & 0 & a_{3,n-\frac{k+3}{2}} & a_{3,n-\frac{k+1}{2}} \\ 
\vdots & \vdots &  & {\mathinner{\mkern2mu\raise1pt\hbox{.}\mkern2mu
\raise4pt\hbox{.}\mkern2mu\raise7pt\hbox{.}\mkern1mu}} & \vdots & \vdots \\ 
0 & 0 &  &  &  &  \\ 
\vdots & \vdots & {\mathinner{\mkern2mu\raise1pt\hbox{.}\mkern2mu
\raise4pt\hbox{.}\mkern2mu\raise7pt\hbox{.}\mkern1mu}} &  &  & {%
\mathinner{\mkern2mu\raise1pt\hbox{.}\mkern2mu
\raise4pt\hbox{.}\mkern2mu\raise7pt\hbox{.}\mkern1mu}} \\ 
0 & a_{n-\frac{k+1}{2},2} &  &  & {\mathinner{\mkern2mu\raise1pt\hbox{.}%
\mkern2mu \raise4pt\hbox{.}\mkern2mu\raise7pt\hbox{.}\mkern1mu}} &  \\ 
a_{n-\frac{k-1}{2},1} & a_{n-\frac{k-1}{2},2} &  &  & {\mathinner{\mkern2mu%
\raise1pt\hbox{.}\mkern2mu
\raise4pt\hbox{.}\mkern2mu\raise7pt\hbox{.}\mkern1mu}} &  \\ 
\vdots & \vdots & {\mathinner{\mkern2mu\raise1pt\hbox{.}\mkern2mu
\raise4pt\hbox{.}\mkern2mu\raise7pt\hbox{.}\mkern1mu}} &  &  &  \\ 
a_{n-1,1} & a_{n-1,2} & \cdots & a_{n-1,\frac{k+1}{2}} & a_{n-1,\frac{k+3}{2}%
} & 0 \\ 
a_{n1} & a_{n2} & \cdots & a_{n,\frac{k+1}{2}} & 0 & 0%
\end{array}%
\right.  \notag \\
&&\ \ \ \ \ \ \ \ \ \ \ \ \ \ \ \ \ \ \left. 
\begin{array}{ccccc}
a_{1,n-\frac{k-1}{2}} & \cdots & a_{1,n-2} & a_{1,n-1} & a_{1,n} \\ 
a_{2,n-\frac{k-1}{2}} & \cdots & a_{2,n-2} & a_{2,n-1} & a_{2,n} \\ 
a_{3,n-\frac{k-1}{2}} & \cdots & a_{3,n-2} & a_{3,n-1} & a_{3,n} \\ 
\vdots & {\mathinner{\mkern2mu\raise1pt\hbox{.}\mkern2mu
\raise4pt\hbox{.}\mkern2mu\raise7pt\hbox{.}\mkern1mu}} &  & \cdots & \vdots
\\ 
&  & a_{\frac{k+1}{2},n-2} & a_{\frac{k+1}{2},n-1} & a_{\frac{k+1}{2},n} \\ 
&  &  &  & 0 \\ 
&  &  & {\mathinner{\mkern2mu\raise1pt\hbox{.}\mkern2mu
\raise4pt\hbox{.}\mkern2mu\raise7pt\hbox{.}\mkern1mu}} & \vdots \\ 
{\mathinner{\mkern2mu\raise1pt\hbox{.}\mkern2mu
\raise4pt\hbox{.}\mkern2mu\raise7pt\hbox{.}\mkern1mu}} &  &  &  & \vdots \\ 
&  &  &  &  \\ 
0 & \cdots & 0 & 0 & 0 \\ 
0 & \cdots & 0 & 0 & a_{nn}%
\end{array}%
\right] .  \label{2}
\end{eqnarray}

The periodic $k$ banded matrices are needed in science and
engineering\linebreak applications for example solving differential
equations by using the finite\linebreak differences, the solution of partial
differential equations and the solution of boundary value problems.

In \cite{1}, the authotrs obtained an algorithm to find the inverse of
the\linebreak periodic tridiagonal matrix by using Doolittle $LU$
factorization and the\linebreak inverse of periodic anti-tridiagonal matrix
is obtained when the inverse of\linebreak periodic tridiagonal matrix
exists. A new algorithm is obtained for the inverse of periodic
pentadiagonal and anti-pentadiagonal matrix in \cite{2}. This paper is an
expansion of \cite{1}. In \cite{3} an algorithm for solving linear systems
having\linebreak periodic pentadiagonal coefficient matrices is obtained. It
is presented that a new computational algorithm to evaluate the determinant
of the tridiagonal matrix with its cost In \cite{4}. In \cite{5} Hadj and
Elouafi obtained a fast numerical algorithm for the inverse of a tridiagonal
and pentadiagonal matrix.

In this work we obtain an algorithm to find the inverse of periodic $k$
banded matrix when its inversion exists. When the algorithm is obtained the
Doolittle $LU$ factorization is used. After finding the inverse of periodic $%
k$ banded matrix, the periodic anti $k$ banded matrix is inverted by using
the inversion of periodic $k$ banded matrix. Also, an algorithm is studied
to solve the linear systems having these coefficient matrices.

\section{Main Result}

In this section, the $LU$ factorization of the matrix $M$ is computed
firstly where $L$ and $U$ are lower and upper triangular matrices,
respectively. It is as in the following:

\begin{equation*}
L=\left[ 
\begin{array}{ccccccccc}
1 &  &  &  &  &  &  &  &  \\ 
l_{21} & \ddots &  &  &  &  &  &  &  \\ 
l_{31} & \ddots &  &  &  &  &  &  &  \\ 
\vdots &  & \ddots &  &  &  &  &  &  \\ 
l_{\frac{k+1}{2},1} & \cdots & l_{\frac{k+1}{2},\frac{k-1}{2}} & 1 &  &  & 
&  &  \\ 
0 & \cdots & l_{\frac{k+3}{2},\frac{k-1}{2}} & l_{\frac{k+3}{2},\frac{k+1}{2}%
} & 1 &  &  &  &  \\ 
0 & \cdots & l_{\frac{k+5}{2},\frac{k-1}{2}} & l_{\frac{k+5}{2},\frac{k+1}{2}%
} & l_{\frac{k+5}{2},\frac{k+3}{2}} & 1 &  &  &  \\ 
\vdots & \ddots & \vdots & \vdots & \vdots & \ddots & \ddots &  &  \\ 
0 & \cdots & 0 & 0 & 0 & \cdots & l_{n-1,n-2} & 1 &  \\ 
l_{n1} & \cdots & l_{n,n-\frac{k+1}{2}} & l_{n,n-\frac{k+3}{2}} & l_{n,n-%
\frac{k+5}{2}} & \cdots & l_{n,n-2} & l_{n,n-1} & 1%
\end{array}%
\right]
\end{equation*}%
and%
\begin{equation*}
U=\left[ 
\begin{array}{ccccccccc}
u_{11} & u_{12} & u_{13} & \cdots & u_{1,\frac{k+1}{2}} & 0 & \cdots & 0 & 
u_{1n} \\ 
& u_{22} & u_{23} & \cdots & u_{2,\frac{k+1}{2}} & u_{2,\frac{k+3}{2}} & 
\cdots & 0 & u_{2n} \\ 
&  & u_{33} & \cdots & u_{3,\frac{k+1}{2}} & u_{3,\frac{k+3}{2}} & \cdots & 0
& u_{3n} \\ 
&  &  & \ddots & \vdots & \vdots &  & \vdots & \vdots \\ 
&  &  &  & u_{\frac{k+1}{2},\frac{k+1}{2}} & u_{\frac{k+1}{2},\frac{k+3}{2}}
& \cdots & 0 & u_{\frac{k+1}{2},n} \\ 
&  &  &  &  & \ddots &  &  & \vdots \\ 
&  &  &  &  &  &  & u_{n-1,n-1} & u_{n-1,n} \\ 
&  &  &  &  &  &  &  & u_{nn}%
\end{array}%
\right]
\end{equation*}%
where%
\begin{equation}
u_{i,r}=a_{i,r}-\dsum\limits_{j=\frac{2r-k+1}{2}}^{i-1}l_{ij}u_{jr}~~,i=1,2,%
\ldots n+i-r-1  \label{3}
\end{equation}%
for $r=i+1,i+2,\ldots i+\frac{k-1}{2}$,%
\begin{equation}
l_{i,r}=\frac{1}{u_{rr}}\left( a_{ir}-\dsum\limits_{j=i-\frac{k-1}{2}%
}^{r-1}l_{ij}u_{jr}\right) ~~,i=1+i-r,\ldots n-1  \label{4}
\end{equation}%
for $r=i-1,i-2,\ldots i-\frac{k-1}{2}$ and

\begin{equation}
u_{ii}=\left\{ 
\begin{array}{ll}
a_{ii}-\dsum\limits_{j=i-\frac{k-1}{2}}^{i-1}l_{ij}u_{ji} & ,i=1,2,\ldots
,n-1 \\ 
a_{nn}-\dsum\limits_{j=1}^{n-1}l_{nj}u_{jn} & ,i=n%
\end{array}%
\right.  \label{5}
\end{equation}%
\begin{equation}
u_{i,n}=\left\{ 
\begin{array}{ll}
a_{1n} & ,i=1 \\ 
-\dsum\limits_{j=i-\frac{k-1}{2}}^{i-1}l_{ij}u_{jn} & ,i=2,3,\ldots ,n-\frac{%
k+1}{2} \\ 
a_{i,n}-\dsum\limits_{j=i-\frac{k-1}{2}}^{i-1}l_{ij}u_{jn} & ,i=n+\frac{1-k}{%
2},\ldots ,n-1%
\end{array}%
\right.  \label{6}
\end{equation}%
\begin{equation}
l_{n,i}=\left\{ 
\begin{array}{ll}
\frac{a_{n1}}{u_{11}} & ,i=1 \\ 
-\frac{1}{u_{ii}}\left( \dsum\limits_{j=i-\frac{k-1}{2}}^{i-1}l_{nj}u_{ji}%
\right) & ,i=2,3,\ldots ,n-\frac{k+1}{2} \\ 
\frac{1}{u_{ii}}\left( a_{ni}-\dsum\limits_{j=i-\frac{k-1}{2}%
}^{i-1}l_{nj}u_{ji}\right) & ,i=n+\frac{1-k}{2},\ldots ,n-1%
\end{array}%
\right.  \label{7}
\end{equation}%
for $i,j\leq 0$ then $l_{ij}=0$ and $u_{ij}=0.$

\bigskip From here, it is clear that%
\begin{equation*}
\det (M)=\dprod\limits_{i=1}^{n}u_{ii}.
\end{equation*}

The inverse matrix is computed as in the following, if the matrix $M$ is
nonsingular:

Let $C_{r}$s be the $r$th column of $M^{-1}$ for $r=1,2,\ldots ,n$ then 
\begin{equation*}
M^{-1}=(S_{i,j})_{1\leq i,j\leq n}\linebreak =(C_{1},C_{2},\ldots
,C_{r},\ldots ,C_{n}).C_{r}=(S_{1,r},S_{2,r},\ldots ,S_{n,r})
\end{equation*}
and it can be written as in the following:

\begin{equation*}
C_{r}=(C_{1},C_{2},\ldots ,C_{r},\ldots ,C_{n})E_{r}
\end{equation*}%
where $E_{r}$ is the Kronecker symbol ($E_{r}=(\delta _{1r,}\delta
_{2r},\ldots ,\delta _{nr,})^{T},r=1,2,\ldots ,n.$).

Now, the algortihm for the inverse of the periodic $k$ banded matrix can be
developed. By using the $LU$ factorization, the entries of the last $\frac{%
k+1}{2}$ columns of $M^{-1}$ are computed as follows:

For $r=n,n-1,\ldots ,m+1,m$ $(m=n-\frac{k-1}{2})$%
\begin{equation}
S_{i,r}=\left\{ 
\begin{array}{cc}
\frac{1}{u_{ii}}\left( t_{ir}-\dsum\limits_{j=i+1}^{n}u_{ij}S_{jr}\right) & 
,i=n,n-1,\ldots ,r+1 \\ 
\frac{1}{u_{ii}}\left( 1-\dsum\limits_{j=i+1}^{n}u_{ij}S_{jr}\right) & ,i=r
\\ 
-\frac{1}{u_{ii}}\left( \dsum\limits_{j=i+1}^{n}u_{ij}S_{jr}\right) & 
,i=r-1,r-2,\ldots ,2,1%
\end{array}%
\right.  \label{8}
\end{equation}%
where $t_{ij}$s are the entries of the inverse of $L$ and for $r=1,2,\ldots
,n-1$ they are computed with the following reccurence relation:%
\begin{equation}
t_{ir}=-l_{ir}-\dsum\limits_{j=r+1}^{i-1}l_{ij}t_{jr}~~,i=r+1,r+2,\ldots ,n.
\label{9}
\end{equation}

Up to now, the entries of the last $\frac{k+1}{2}$ columns are obtained.
Entries of the remaining $n-\frac{k-1}{2}$ columns are computed by using the
following equation%
\begin{equation*}
M^{-1}M=I_{n}.
\end{equation*}

For $j=n-\frac{k+1}{2},n-\frac{k+3}{2},\ldots ,1$%
\begin{equation}
C_{j}=\frac{1}{a_{j,j+\frac{k-1}{2}}}\left( E_{j+\frac{k-1}{2}%
}-\dsum\limits_{r=j+1}^{j+k-1}a_{r,j+\frac{k-1}{2}}C_{r}\right)  \label{10}
\end{equation}%
where if $j=1,2,\ldots ,n-\frac{k+1}{2}$ then $a_{j,j+\frac{k-1}{2}}\neq 0.$
Here, $i>n,\ j>n$ then $a_{ij}=0,$ if $i>n$ then $C_{i}=0.$

\textbf{Algorithm 1:}

\textbf{INPUT:} $n$ is the order of the periodic $k$ banded matrix, $k$ is
the bandwidth of the matrix.

\textbf{OUTPUT: }The inverse matrix $M^{-1}=\left( S_{i,j}\right) _{1\leq
i,j\leq n}.$

\textbf{Step1:} For $i=1,2,\ldots ,n-\frac{k+1}{2}$, if $a_{i,i+\frac{k-1}{2}%
}=0,$ then $a_{i,i+\frac{k-1}{2}}=\lambda .$

\textbf{Step2: }For $i=\frac{k+3}{2},\frac{k+5}{2},\ldots ,n$, if $a_{i,i-%
\frac{k-1}{2}}=0,$ then $a_{i,i-\frac{k-1}{2}}=\lambda .$

\textbf{Step3: }By using (\ref{3})-(\ref{7}), compute the elements of the
matrices $L$ and $U.$ For $i=1,2,\ldots ,n,$ if $u_{ii}=0,$ then $%
u_{ii}=\lambda .$

\textbf{Step 4:} Compute $\det (M)=\left(
\dsum\limits_{i=1}^{n}u_{ii}\right) _{\lambda =0}.$ If $M$ is singular, then
the output is "Singular Matrix".

\textbf{Step 5:} Compute the elements $t_{i,r}$ by using (\ref{9}).

\textbf{Step6:} Compute the elements of the last $\frac{k+1}{2}$ columns by
using (\ref{8}).

\textbf{Step7:} Compute the elements of the remaining $n-\frac{k+1}{2}$
columns by using (\ref{10}).

\textbf{Step8:} Substitute the actual value of $\lambda $ in all elements of
the inverse matrix $M^{-1}.$

The inverse matrix of the periodic anti $k$ banded matrix $N$ can be
obtained by using the inverse of periodic $k$ banded matrix $M$ easily.

Let $R$ be an $n\times n$ matrix as in the following form:

\begin{equation*}
R=\left[ 
\begin{array}{ccccc}
0 & \cdots & \cdots & 0 & 1 \\ 
\vdots &  & {\mathinner{\mkern2mu\raise1pt\hbox{.}\mkern2mu
\raise4pt\hbox{.}\mkern2mu\raise7pt\hbox{.}\mkern1mu}} & 1 & 0 \\ 
\vdots & {\mathinner{\mkern2mu\raise1pt\hbox{.}\mkern2mu
\raise4pt\hbox{.}\mkern2mu\raise7pt\hbox{.}\mkern1mu}} & {%
\mathinner{\mkern2mu\raise1pt\hbox{.}\mkern2mu
\raise4pt\hbox{.}\mkern2mu\raise7pt\hbox{.}\mkern1mu}} & {%
\mathinner{\mkern2mu\raise1pt\hbox{.}\mkern2mu
\raise4pt\hbox{.}\mkern2mu\raise7pt\hbox{.}\mkern1mu}} & \vdots \\ 
0 & {\mathinner{\mkern2mu\raise1pt\hbox{.}\mkern2mu
\raise4pt\hbox{.}\mkern2mu\raise7pt\hbox{.}\mkern1mu}} & {%
\mathinner{\mkern2mu\raise1pt\hbox{.}\mkern2mu
\raise4pt\hbox{.}\mkern2mu\raise7pt\hbox{.}\mkern1mu}} &  & \vdots \\ 
1 & 0 & \cdots & \cdots & 0%
\end{array}%
\right]
\end{equation*}%
It is clear that $R$ is nonsingular and its inversion is itself \cite{2}.

The following relation is true for (\ref{1})and (\ref{2}).

\begin{equation*}
N=MR
\end{equation*}%
Thus, the inverse of (\ref{2}) is obtained as in the following:%
\begin{equation*}
N^{-1}=RM^{-1}.
\end{equation*}

Also, linear systems having periodic $k$ banded coefficient matrix, $Mx=y,$
can be solved by using the Dolittle $LU$ factorization of this type of
matrix. Here the component $u_{ii}$ is important. Because, the system has a
unique solution if $u_{ii}\neq 0,\ \ \ i=1,2,\ldots ,n.$

\textbf{Algorithm 2:}

\textbf{Step1: }By using the Step1-Step3 of Algorithm1, compute the $LU$%
\linebreak factorization of the coefficient matrix $M.$

\textbf{Step2: }For $i=1,2,\ldots ,n$ \ \ $z_{i}=y_{i}-\dsum%
\limits_{j=1}^{i-1}l_{ij}z_{j}.$

\textbf{Step3:} For $i=n,n-1,\ldots ,2,1$%
\begin{equation*}
x_{i}=\left\{ 
\begin{array}{ll}
\frac{1}{u_{ii}}z_{i} & ,i=n \\ 
\frac{1}{u_{ii}}\left( z_{i}-\dsum\limits_{j=i+1}^{i+\frac{k-1}{2}%
}u_{ij}x_{j}\right) & ,i=n-1,\ldots ,n-\frac{k-1}{2} \\ 
\frac{1}{u_{ii}}\left( z_{i}-\dsum\limits_{j=i+1}^{i+\frac{k-1}{2}%
}u_{ij}x_{j}-u_{in}x_{n}\right) & ,i=n-\frac{k+1}{2},\ldots ,1%
\end{array}%
\right. .
\end{equation*}

\section{Numerical Example}

\begin{example}
Consider the $6\times 6$ matrices $M$ and $N$ as in the following%
\begin{equation*}
M=\left[ 
\begin{array}{rrrrrr}
2 & 1 & 0 & 0 & 0 & 1 \\ 
1 & -1 & 2 & 0 & 0 & 0 \\ 
0 & 2 & -2 & 3 & 0 & 0 \\ 
0 & 0 & -1 & 1 & 1 & 0 \\ 
0 & 0 & 0 & 2 & -3 & -2 \\ 
2 & 0 & 0 & 0 & 1 & 5%
\end{array}%
\right]
\end{equation*}%
and%
\begin{equation*}
N=\left[ 
\begin{array}{rrrrrr}
1 & 0 & 0 & 0 & 1 & 2 \\ 
0 & 0 & 0 & 2 & -1 & 1 \\ 
0 & 0 & 3 & -2 & 2 & 0 \\ 
0 & 1 & 1 & -1 & 0 & 0 \\ 
-2 & -3 & 2 & 0 & 0 & 0 \\ 
5 & 1 & 0 & 0 & 0 & 2%
\end{array}%
\right] .
\end{equation*}%
We apply the Algorithm1 to the matrix $M$ and we have
\end{example}

\begin{itemize}
\item For $i=1,2,3,4$ $\ \ \ \
u_{i,i+1}=(u_{12},u_{23},u_{34},u_{45})=(1,2,3,1)$

\item For $i=2,3,4,5\ \ \ \ \ l_{i,i-1}=(l_{21},l_{32},l_{43},l_{54})=(\frac{%
1}{2},-\frac{4}{3},-\frac{3}{2},\frac{4}{11})$

\item For $i=1,2,3,4,5,6\ \ \ \ \
u_{ii}=(u_{11},u_{22},u_{33},u_{44},u_{55},u_{66})=(2,-\frac{3}{2},\frac{2}{3%
},\frac{11}{2},-\frac{37}{11},\frac{153}{37})$

\item For $i=1,2,3,4,5\ \ \ \ \
u_{i,6}=(u_{16},u_{26},u_{36},u_{46},u_{56})=(1,-\frac{1}{2},-\frac{2}{3}%
,-1,-\frac{18}{11})$

\item For$\ i=1,2,3,4,5\ \ \ \ \
l_{6,i}=(l_{61},l_{62},l_{63},l_{64},l_{65})=(1,\frac{2}{3},-2,\frac{12}{11},%
\frac{1}{37})$

\item $\det (M)=153$

\item For $i=2,3,4,5,6\ \ \ \ \
t_{i,1}=(t_{21},t_{31},t_{41},t_{51},t_{61})=(-\frac{1}{2},-\frac{2}{3},-1,%
\frac{4}{11},-\frac{34}{37})$

\item For $i=3,4,5,6\ \ \ \ \ t_{i,2}=(t_{32},t_{42},t_{52},t_{62})=(\frac{4%
}{3},2,-\frac{8}{11},-\frac{6}{37})$

\item For $i=4,5,6\ \ \ \ \ t_{i,3}=(t_{43},t_{53},t_{63})=(\frac{3}{2},-%
\frac{6}{11},\frac{14}{37})$

\item For $i=5,6\ \ \ \ \ t_{i,4}=(t_{54},t_{64})=(-\frac{4}{11},-\frac{40}{%
37})$

\item For $i=6\ \ \ \ \ t_{i,5}=(t_{65})=(-\frac{1}{37})$

\item $C_{6}=(S_{66},S_{56},S_{46},S_{36},S_{26},S_{16})=(\frac{37}{153},-%
\frac{2}{17},\frac{10}{153},-\frac{8}{153},-\frac{23}{153},-\frac{7}{153})$

\item $C_{5}=(S_{65},S_{55},S_{45},S_{35},S_{25},S_{15})=(-\frac{1}{153},-%
\frac{5}{17},\frac{8}{153},-\frac{37}{153},-\frac{49}{153},\frac{25}{153})$

\item $C_{4}=(S_{64},S_{54},S_{44},S_{34},S_{24},S_{14})=(-\frac{40}{153},%
\frac{4}{17},\frac{14}{153},-\frac{103}{153},-\frac{124}{153},\frac{82}{153}%
) $

\item $C_{3}=(S_{63},S_{53},S_{43},S_{33},S_{23},S_{13})=(\frac{14}{153},%
\frac{2}{17},\frac{41}{153},\frac{59}{153},\frac{74}{153},-\frac{44}{153})$

\item $C_{2}=(S_{62},S_{52},S_{42},S_{32},S_{22},S_{12})=(-\frac{2}{51},%
\frac{4}{17},\frac{16}{51},\frac{28}{51},\frac{4}{51},-\frac{1}{51})$

\item $C_{1}=(S_{61},S_{51},S_{41},S_{31},S_{21},S_{11})=(-\frac{2}{9},0,-%
\frac{2}{9},-\frac{2}{9},\frac{1}{9},\frac{5}{9})$

\item $M^{-1}=\left[ 
\begin{array}{rrrrrr}
\frac{5}{9} & -\frac{1}{51} & -\frac{144}{153} & \frac{82}{153} & \frac{25}{%
153} & -\frac{7}{153} \\ 
\frac{1}{9} & \frac{4}{51} & \frac{74}{153} & -\frac{124}{153} & -\frac{49}{%
153} & -\frac{23}{153} \\ 
-\frac{2}{9} & \frac{28}{51} & \frac{59}{153} & -\frac{103}{153} & -\frac{37%
}{153} & -\frac{8}{153} \\ 
-\frac{2}{9} & \frac{16}{51} & \frac{41}{153} & \frac{14}{153} & \frac{8}{153%
} & \frac{10}{153} \\ 
0 & \frac{4}{17} & \frac{2}{17} & \frac{4}{17} & -\frac{5}{17} & -\frac{2}{17%
} \\ 
-\frac{2}{9} & -\frac{2}{51} & \frac{14}{153} & -\frac{40}{153} & -\frac{1}{%
153} & \frac{37}{153}%
\end{array}%
\right] $

\item $N^{-1}=\left[ 
\begin{array}{cccccc}
0 & 0 & 0 & 0 & 0 & 1 \\ 
0 & 0 & 0 & 0 & 1 & 0 \\ 
0 & 0 & 0 & 1 & 0 & 0 \\ 
0 & 0 & 1 & 0 & 0 & 0 \\ 
0 & 1 & 0 & 0 & 0 & 0 \\ 
1 & 0 & 0 & 0 & 0 & 0%
\end{array}%
\right] \left[ 
\begin{array}{rrrrrr}
\frac{5}{9} & -\frac{1}{51} & -\frac{144}{153} & \frac{82}{153} & \frac{25}{%
153} & -\frac{7}{153} \\ 
\frac{1}{9} & \frac{4}{51} & \frac{74}{153} & -\frac{124}{153} & -\frac{49}{%
153} & -\frac{23}{153} \\ 
-\frac{2}{9} & \frac{28}{51} & \frac{59}{153} & -\frac{103}{153} & -\frac{37%
}{153} & -\frac{8}{153} \\ 
-\frac{2}{9} & \frac{16}{51} & \frac{41}{153} & \frac{14}{153} & \frac{8}{153%
} & \frac{10}{153} \\ 
0 & \frac{4}{17} & \frac{2}{17} & \frac{4}{17} & -\frac{5}{17} & -\frac{2}{17%
} \\ 
-\frac{2}{9} & -\frac{2}{51} & \frac{14}{153} & -\frac{40}{153} & -\frac{1}{%
153} & \frac{37}{153}%
\end{array}%
\right] $

$\ \ \ \ \ \ =\left[ 
\begin{array}{rrrrrr}
-\frac{2}{9} & -\frac{2}{51} & \frac{14}{153} & -\frac{40}{153} & -\frac{1}{%
153} & \frac{37}{153} \\ 
0 & \frac{4}{17} & \frac{2}{17} & \frac{4}{17} & -\frac{5}{17} & -\frac{2}{17%
} \\ 
-\frac{2}{9} & \frac{16}{51} & \frac{41}{153} & \frac{14}{153} & \frac{8}{153%
} & \frac{10}{153} \\ 
-\frac{2}{9} & \frac{28}{51} & \frac{59}{153} & -\frac{103}{153} & -\frac{37%
}{153} & -\frac{8}{153} \\ 
\frac{1}{9} & \frac{4}{51} & \frac{74}{153} & -\frac{124}{153} & -\frac{49}{%
153} & -\frac{23}{153} \\ 
\frac{5}{9} & -\frac{1}{51} & -\frac{16}{17} & \frac{82}{153} & \frac{25}{153%
} & -\frac{7}{153}%
\end{array}%
\right] \allowbreak $
\end{itemize}

\begin{example}
Consider the $10\times 10$ matrices $M$ and $N$ as in the following%
\begin{equation*}
M=\left[ 
\begin{array}{rrrrrrrrrr}
1 & -1 & 2 & 2 & -1 & 0 & 0 & 0 & 0 & 1 \\ 
2 & -1 & 3 & 1 & 1 & 2 & 0 & 0 & 0 & 0 \\ 
1 & -1 & 1 & 2 & 1 & -2 & -1 & 0 & 0 & 0 \\ 
-3 & 1 & -1 & 1 & -3 & 1 & 1 & -3 & 0 & 0 \\ 
2 & -1 & 1 & 0 & -3 & 2 & 1 & -1 & -1 & 0 \\ 
0 & 1 & 2 & 0 & -1 & 0 & -2 & 1 & 0 & 1 \\ 
0 & 0 & -2 & 0 & 1 & -1 & 1 & -2 & 1 & -1 \\ 
0 & 0 & 0 & 1 & 3 & 2 & -1 & 1 & 2 & 1 \\ 
0 & 0 & 0 & 0 & -1 & 0 & 2 & 1 & -2 & 1 \\ 
2 & 0 & 0 & 0 & 0 & 2 & 1 & 1 & -1 & 2%
\end{array}%
\right]
\end{equation*}%
and%
\begin{equation*}
N=\left[ 
\begin{array}{rrrrrrrrrr}
2 & 0 & 0 & 0 & 0 & 2 & 1 & 1 & -1 & 2 \\ 
0 & 0 & 0 & 0 & -1 & 0 & 2 & 1 & -2 & 1 \\ 
0 & 0 & 0 & 1 & 3 & 2 & -1 & 1 & 2 & 1 \\ 
0 & 0 & -2 & 0 & 1 & -1 & 1 & -2 & 1 & -1 \\ 
0 & 1 & 2 & 0 & -1 & 0 & -2 & 1 & 0 & 1 \\ 
2 & -1 & 1 & 0 & -3 & 2 & 1 & -1 & -1 & 0 \\ 
-3 & 1 & -1 & 1 & -3 & 1 & 1 & -3 & 0 & 0 \\ 
1 & -1 & 1 & 2 & 1 & -2 & -1 & 0 & 0 & 0 \\ 
2 & -1 & 3 & 1 & 1 & 2 & 0 & 0 & 0 & 0 \\ 
1 & -1 & 2 & 2 & -1 & 0 & 0 & 0 & 0 & 1%
\end{array}%
\right]
\end{equation*}%
We apply the Algorithm1 to the matrix $M$ and we have
\end{example}

\begin{itemize}
\item For $i=1,2,3,4,5,6,7,8\ \ \ \ \ \ \ \ \ \ \ \ \ \ \ \ \ \ \ \ \ \ \ \
\ \ \ \ \ \ \ \ \ \ \ \ \ \ \ \ \ \ \ \ \ \ \ \ \ \ \ \ \ \ \ \ \ \ \ \ \ \
\ \ \ \ \ \ \ \ \ \ \ \ \ \ \ \ \ \ \ \ \ \ \ \ \ \ \ \ \ \ \ \ \ \ \ \ \ \
\ \ $\linebreak $%
u_{i,i+1}=(u_{12},u_{23},u_{34},u_{45},u_{56,}u_{67},u_{78},u_{89})=(-1,-1,0,6,3,-7,-%
\frac{53}{29},\frac{175}{54})$

\item For $i=1,2,3,4,5,6,7\ \ \ \ \ \ \ \ \ \ \ \ \ \ \ \ \ \ \ \ \ \ \ \ \
\ \ \ \ \ \ \ \ \ \ \ \ \ \ \ \ \ \ \ \ \ \ \ \ \ \ \ \ \ \ \ \ \ \ \ \ \ \
\ \ \ \ \ \ \ \ \ \ \ \ \ \ \ \ \ \ \ \ \ \ \ \ \ \ \ \ $\linebreak $%
u_{i,i+2}=(u_{13},u_{24},u_{35},u_{46},u_{57},u_{68},u_{79})=(2,-3,2,-1,1,42,%
\frac{121}{58})$

\item For $i=1,2,3,4,5,6\ \ \ \ \
u_{i,i+3}=(u_{14},u_{25},u_{36},u_{47},u_{58},u_{69})=(2,3,-2,-2,-4,8)$

\item For $i=1,2,3,4,5\ \ \ \ \
u_{i,i+4}=(u_{15},u_{26},u_{37},u_{48},u_{59})=(-1,2,-1,-3,-1)$

\item For $i=2,3,4,5,6,7,8,9\ \ \ \ \ \ \ \ \ \ \ \ \ \ \ \ \ \ \ \ \ \ \ \
\ \ \ \ \ \ \ \ \ \ \ \ \ \ \ \ \ \ \ \ \ \ \ \ \ \ \ \ \ \ \ \ \ \ \ \ \ \
\ \ \ \ \ \ \ \ \ \ \ \ \ \ \ \ \ \ \ \ \ \ \ \ \ \ \ \ \ \ \ \ \ \ \ \ \ \
\ $\linebreak $%
l_{i,i-1}=(l_{21},l_{32},l_{43},l_{54},l_{65},l_{76},l_{87},l_{98})=(2,0,-3,-1,8,%
\frac{3}{58},-\frac{2}{27},-\frac{27}{205})$

\item For $i=3,4,5,6,7,8,9\ \ \ \ \
l_{i,i-2}=(l_{31},l_{42},l_{53},l_{64},l_{75},l_{86},l_{97})=(1,-2,2,3,\frac{%
3}{2},\frac{3}{58},1)$

\item For $i=4,5,6,7,8,9\ \ \ \ \
l_{i,i-3}=(l_{41},l_{52},l_{63},l_{74},l_{85},l_{96})=(-3,1,-3,0,\frac{3}{2},%
\frac{3}{58})$

\item For $i=5,6,7,8,9\ \ \ \ \
l_{i,i-4}=(l_{51},l_{62},l_{73},l_{84},l_{95})=(2,1,2,1,\frac{1}{2})$

\item For $i=1,2,3,4,5,6,7,8,9,10\ \ \ \ \ \ \ \ \ \ \ \ \ \ \ \ \ \ \ \ \ \
\ \ \ \ \ \ \ \ \ \ \ \ \ \ \ \ \ \ \ \ \ \ \ \ \ \ \ \ \ \ \ \ \ \ \ \ \ \
\ \ \ \ \ \ \ \ \ \ \ \ \ \ \ \ \ \ \ \ \ \ \ \ \ \ \ \ \ \ \ \ \ \ \ \ \ \
\ \ $\linebreak $\
u_{ii}=(u_{11},u_{22},u_{33},u_{44},u_{55},u_{66},u_{77},u_{88},u_{99},u_{10,10}) 
$

$\ \ \ \ \ =(1,1,-1,1,-2,-29,\frac{54}{29},\frac{215}{27},-\frac{309}{86},%
\frac{944}{1545})$

\item For $i=1,2,3,4,5,6,7,8,9\ \ \ \ \ \ \ \ \ \ \ \ \ \ \ \ \ \ \ \ \ \ \
\ \ \ \ \ \ \ \ \ \ \ \ \ \ \ \ \ \ \ \ \ \ \ \ \ \ \ \ \ \ \ \ \ \ \ \ \ \
\ \ \ \ \ \ \ \ \ \ \ \ \ \ \ \ \ \ \ \ \ \ \ \ \ \ \ \ \ \ \ \ \ \ \ \ \ \
\ \ \ $\linebreak $%
u_{i,10}=(u_{1,10},u_{2,10},u_{3,10},u_{4,10},u_{5,10},u_{6,10},u_{7,10},u_{8,10},u_{9,10},) 
$

$\ \ \ \ \ \ \ =(1,-2,-1,-4,-2,28,\frac{74}{29},\frac{182}{27},-\frac{248}{%
215})$

\item For$\ i=1,2,3,4,5,6,7,8,9\ \ \ \ \ \ \ \ \ \ \ \ \ \ \ \ \ \ \ \ \ \ \
\ \ \ \ \ \ \ \ \ \ \ \ \ \ \ \ \ \ \ \ \ \ \ \ \ \ \ \ \ \ \ \ \ \ \ \ \ \
\ \ \ \ \ \ \ \ \ \ \ \ \ \ \ \ \ \ \ \ \ \ \ \ \ \ \ \ \ \ \ \ \ \ \ \ \ \
\ \ \ \ \ \ $\linebreak $%
l_{10,i}=(l_{10,1},l_{10,2},l_{10,3},l_{10,4},l_{10,5},l_{10,6},l_{10,7},l_{10,8},l_{10,9}) 
$

$\ \ \ \ \ \ =(2,2,2,2,10,\frac{26}{29},\frac{95}{54},\frac{331}{430},\frac{%
373}{309})$

\item $\det (M)=1888$

\item For $i=2,3,4,5,6,7,8,9,10\ \ \ \ \ \ \ \ \ \ \ \ \ \ \ \ \ \ \ \ \ \ \
\ \ \ \ \ \ \ \ \ \ \ \ \ \ \ \ \ \ \ \ \ \ \ \ \ \ \ \ \ \ \ \ \ \ \ \ \ \
\ \ \ \ \ \ \ \ \ \ \ \ \ \ \ \ \ \ \ \ \ \ \ \ \ \ \ \ \ \ \ \ \ \ \ \ \ \
\ \ $\linebreak $%
t_{i,1}=(t_{21},t_{31},t_{41},t_{51},t_{61},t_{71},t_{81},t_{91},t_{10,1})$

$\ \ \ \ \ =(-2,-1,-4,-2,27,\frac{209}{58},\frac{317}{54},-\frac{1403}{430},%
\frac{2699}{3090})$

\item For $i=3,4,5,6,7,8,9,10\ \ \ \ \ \ \ \ \ \ \ \ \ \ \ \ \ \ \ \ \ \ \ \
\ \ \ \ \ \ \ \ \ \ \ \ \ \ \ \ \ \ \ \ \ \ \ \ \ \ \ \ \ \ \ \ \ \ \ \ \ \
\ \ \ \ \ \ \ \ \ \ \ \ \ \ \ \ \ \ \ \ \ \ \ \ \ \ \ \ \ \ \ \ \ \ \ \ \ \
\ \ \ \ \ \ \ \ \ \ \ \ \ \ \ \ \ \ \ \ \ \ \ \ \ $\linebreak $%
t_{i,2}=(t_{32},t_{42},t_{52},t_{62},t_{72},t_{82},t_{92},t_{10,2})=(0,2,1,-15,-%
\frac{21}{29},-\frac{25}{9},\frac{28}{43},\frac{23}{209})$

\item For $i=4,5,6,7,8,9,10$

$t_{i,3}=(t_{43},t_{53},t_{63},t_{73},t_{83},t_{93},t_{10,3})=(3,1,-14,-%
\frac{161}{58},-\frac{215}{54},\frac{5}{2},-\frac{160}{309})$

\item For $i=5,6,7,8,9,10$

$t_{i,4}=(t_{54},t_{64},t_{74},t_{84},t_{94},t_{10,4})=(1,-11,-\frac{27}{29}%
,-2,\frac{161}{215},\frac{419}{3090})$

\item For $i=6,7,8,9,10$

$t_{i,5}=(t_{65},t_{75},t_{85},t_{95},t_{10,5})=(-8,-\frac{63}{58},-\frac{7}{%
6},\frac{367}{430},-\frac{3241}{3090})$

\item For $i=7,8,9,10$

$t_{i,6}=(t_{76},t_{86},t_{96},t_{10,6})=(-\frac{3}{58},-\frac{1}{18},-\frac{%
3}{430},\frac{777}{1030})$

\item For $i=8,9,10\ \ \ \ \ t_{i,7}=(t_{87},t_{97},t_{10,7})=(\frac{2}{27},-%
\frac{213}{215},-\frac{639}{1030})$

\item For $i=9,10\ \ \ \ \ t_{i,8}=(t_{98},t_{10,8})=(\frac{27}{215},-\frac{%
949}{1030})$

\item For $i=10\ \ \ \ \ t_{i,9}=(t_{10,9})=(-\frac{373}{309})$

\item $%
C_{10}=(S_{10,10},S_{9,10},S_{8,10},S_{7,10},S_{6,10},S_{5,10},S_{4,10},S_{3,10},S_{2,10},S_{1,10})\ \ \ \ \ \ \ \ \ \ \ \ \ \ \ \ \ \ \ \ \ \ \ \ \ \ \ \ \ \ \ \ \ \ \ \ \ \ \ \ \ \ \ \ \ \ \ \ \ \ \ \ \ \ 
$\linebreak $~~~~~~=(\frac{1545}{944},-\frac{31}{59},-\frac{553}{472},-\frac{%
119}{236},-\frac{33}{236},\frac{479}{944},-\frac{137}{118},\frac{153}{944},-%
\frac{609}{472},-\frac{199}{472})$

\item $%
C_{9}=(S_{10,9},S_{99},S_{89},S_{79},S_{69},S_{59},S_{49},S_{39},S_{29},S_{19})\ \ \ \ \ \ \ \ \ \ \ \ \ \ \ \ \ \ \ \ \ \ \ \ \ \ \ \ \ \ \ \ \ \ \ \ \ \ \ \ \ \ \ \ \ \ \ \ \ \ \ \ \ \ 
$\linebreak $~~~~~=(-\frac{1865}{944},\frac{21}{59},\frac{721}{472},\frac{191%
}{236},\frac{49}{236},-\frac{511}{944},\frac{207}{118},-\frac{313}{944},%
\frac{1033}{472},\frac{367}{472})$

\item $%
C_{8}=(S_{10,8},S_{98},S_{88},S_{78},S_{68},S_{58},S_{48},S_{38},S_{28},S_{18})\ \ \ \ \ \ \ \ \ \ \ \ \ \ \ \ \ \ \ \ \ \ \ \ \ \ \ \ \ \ \ \ \ \ \ \ \ \ \ \ \ \ \ \ \ \ \ \ \ \ \ \ \ \ 
$\linebreak $~~~~~=(-\frac{2847}{1888},\frac{53}{118},\frac{1151}{944},\frac{%
173}{472},\frac{163}{472},-\frac{857}{1888},\frac{337}{286},-\frac{863}{1888}%
,\frac{1399}{944},\frac{561}{944})$

\item $%
C_{7}=(S_{10,7},S_{97},S_{87},S_{77},S_{67},S_{57},S_{47},S_{37},S_{27},S_{17})\ \ \ \ \ \ \ \ \ \ \ \ \ \ \ \ \ \ \ \ \ \ \ \ \ \ \ \ \ \ \ \ \ \ \ \ \ \ \ \ \ \ \ \ \ \ \ \ \ \ \ \ \ \ 
$\linebreak $~~~~~=(-\frac{1917}{1888},\frac{71}{118},\frac{589}{944},\frac{%
303}{472},-\frac{31}{472},-\frac{587}{1888},\frac{211}{236},-\frac{221}{1888}%
,\frac{1509}{944},\frac{707}{944})$

\item $%
C_{6}=(S_{10,6},S_{96},S_{86},S_{76},S_{66},S_{56},S_{46},S_{36},S_{26},S_{16})\ \ \ \ \ \ \ \ \ \ \ \ \ \ \ \ \ \ \ \ \ \ \ \ \ \ \ \ \ \ \ \ \ \ \ \ \ \ \ \ \ \ \ \ \ \ \ \ \ \ \ \ \ \ 
$\linebreak $~~~~~=(-\frac{2331}{1888},\frac{47}{118},\frac{827}{944},\frac{%
169}{472},\frac{31}{472},-\frac{829}{1888},\frac{261}{236},-\frac{251}{1888},%
\frac{1795}{944},\frac{709}{944})$

\item $%
C_{5}=(S_{10,5},S_{95},S_{85},S_{75},S_{65},S_{55},S_{45},S_{35},S_{25},S_{15})\ \ \ \ \ \ \ \ \ \ \ \ \ \ \ \ \ \ \ \ \ \ \ \ \ \ \ \ \ \ \ \ \ \ \ \ \ \ \ \ \ \ \ \ \ \ \ \ \ \ \ \ \ \ 
$\linebreak $~~~~~=(-\frac{3241}{1888},\frac{37}{118},\frac{1113}{944},\frac{%
123}{472},\frac{165}{472},-\frac{1215}{1885},\frac{331}{236},-\frac{1001}{%
1888},\frac{1393}{944},\frac{759}{944})$

\item $%
C_{4}=(S_{10,4},S_{94},S_{84},S_{74},S_{64},S_{54},S_{44},S_{34},S_{24},S_{14})\ \ \ \ \ \ \ \ \ \ \ \ \ \ \ \ \ \ \ \ \ \ \ \ \ \ \ \ \ \ \ \ \ \ \ \ \ \ \ \ \ \ \ \ \ \ \ \ \ \ \ \ \ \ 
$\linebreak $~~~~~=(\frac{419}{1888},-\frac{33}{118},-\frac{307}{944},-\frac{%
81}{472},\frac{41}{472},\frac{213}{1888},-\frac{5}{236},\frac{3}{1888},-%
\frac{123}{944},-\frac{189}{944})$

\item $%
C_{3}=(S_{10,3},S_{93},S_{83},S_{73},S_{63},S_{53},S_{43},S_{33},S_{23},S_{13})\ \ \ \ \ \ \ \ \ \ \ \ \ \ \ \ \ \ \ \ \ \ \ \ \ \ \ \ \ \ \ \ \ \ \ \ \ \ \ \ \ \ \ \ \ \ \ \ \ \ \ \ \ \ 
$\linebreak $~~~~~=(-\frac{50}{59},-\frac{25}{59},\frac{23}{59},-\frac{14}{59%
},\frac{10}{59},\frac{5}{59},\frac{58}{59},-\frac{25}{59},\frac{44}{59},%
\frac{23}{59})$

\item $%
C_{2}=(S_{10,2},S_{92},S_{82},S_{72},S_{62},S_{52},S_{42},S_{32},S_{22},S_{12})\ \ \ \ \ \ \ \ \ \ \ \ \ \ \ \ \ \ \ \ \ \ \ \ \ \ \ \ \ \ \ \ \ \ \ \ \ \ \ \ \ \ \ \ \ \ \ \ \ \ \ \ \ \ 
$\linebreak $~~~~~=(\frac{115}{944},-\frac{13}{59},-\frac{171}{472},\frac{11%
}{236},\frac{9}{236},\frac{277}{944},-\frac{27}{118},\frac{323}{944},-\frac{%
27}{472},-\frac{53}{472})$

\item $%
C_{1}=(S_{10,1},S_{91},S_{81},S_{71},S_{61},S_{51},S_{41},S_{31},S_{21},S_{11})\ \ \ \ \ \ \ \ \ \ \ \ \ \ \ \ \ \ \ \ \ \ \ \ \ \ \ \ \ \ \ \ \ \ \ \ \ \ \ \ \ \ \ \ \ \ \ \ \ \ \ \ \ \ 
$\linebreak $~~~~~=(\frac{2699}{1888},\frac{53}{118},-\frac{619}{944},\frac{%
55}{472},-\frac{191}{472},\frac{205}{1888},-\frac{253}{236},\frac{907}{1888}%
,-\frac{1315}{944},-\frac{501}{944})$

\item $M^{-1}=\left[ 
\begin{array}{rrrrrrrrrr}
-\frac{501}{944} & \frac{53}{472} & \frac{23}{59} & -\frac{189}{944} & \frac{%
759}{944} & \frac{709}{944} & \frac{707}{944} & \frac{561}{944} & \frac{367}{%
472} & -\frac{199}{472} \\ 
-\frac{1315}{944} & -\frac{27}{472} & \frac{44}{59} & -\frac{123}{944} & 
\frac{1393}{944} & \frac{1795}{944} & \frac{1509}{944} & \frac{1399}{944} & 
\frac{1033}{472} & -\frac{609}{472} \\ 
\frac{907}{1888} & \frac{323}{944} & -\frac{25}{59} & \frac{3}{1888} & -%
\frac{1001}{1888} & -\frac{251}{1888} & -\frac{221}{1888} & -\frac{863}{1888}
& -\frac{313}{944} & \frac{153}{944} \\ 
-\frac{253}{236} & -\frac{27}{118} & \frac{58}{59} & -\frac{5}{236} & \frac{%
331}{236} & \frac{261}{236} & \frac{211}{236} & \frac{337}{236} & \frac{207}{%
118} & -\frac{237}{118} \\ 
\frac{205}{118} & \frac{277}{944} & -\frac{5}{59} & \frac{213}{1888} & -%
\frac{1215}{1888} & -\frac{829}{1888} & -\frac{587}{1888} & -\frac{857}{1888}
& -\frac{511}{944} & \frac{479}{944} \\ 
-\frac{191}{472} & \frac{9}{236} & \frac{10}{59} & \frac{41}{472} & \frac{165%
}{472} & \frac{31}{472} & -\frac{31}{472} & \frac{163}{472} & \frac{49}{236}
& -\frac{33}{236} \\ 
\frac{55}{472} & \frac{11}{236} & -\frac{14}{59} & -\frac{81}{472} & \frac{%
123}{472} & \frac{169}{472} & \frac{303}{472} & \frac{173}{472} & \frac{191}{%
236} & -\frac{119}{236} \\ 
-\frac{619}{944} & -\frac{171}{472} & \frac{23}{59} & -\frac{307}{944} & 
\frac{1113}{944} & \frac{827}{944} & \frac{589}{944} & \frac{1151}{944} & 
\frac{721}{472} & -\frac{553}{472} \\ 
\frac{53}{118} & -\frac{13}{59} & -\frac{25}{59} & -\frac{33}{118} & \frac{37%
}{118} & \frac{47}{118} & \frac{71}{118} & \frac{53}{118} & \frac{21}{59} & -%
\frac{31}{59} \\ 
\frac{2699}{1888} & \frac{115}{944} & -\frac{50}{59} & \frac{419}{1888} & -%
\frac{3241}{1888} & -\frac{2331}{1888} & -\frac{1917}{1888} & -\frac{2847}{%
1888} & -\frac{1865}{944} & \frac{1545}{944}%
\end{array}%
\right] $

\item $N^{-1}=\left[ 
\begin{array}{cccccccccc}
0 & 0 & 0 & 0 & 0 & 0 & 0 & 0 & 0 & 1 \\ 
0 & 0 & 0 & 0 & 0 & 0 & 0 & 0 & 1 & 0 \\ 
0 & 0 & 0 & 0 & 0 & 0 & 0 & 1 & 0 & 0 \\ 
0 & 0 & 0 & 0 & 0 & 0 & 1 & 0 & 0 & 0 \\ 
0 & 0 & 0 & 0 & 0 & 1 & 0 & 0 & 0 & 0 \\ 
0 & 0 & 0 & 0 & 1 & 0 & 0 & 0 & 0 & 0 \\ 
0 & 0 & 0 & 1 & 0 & 0 & 0 & 0 & 0 & 0 \\ 
0 & 0 & 1 & 0 & 0 & 0 & 0 & 0 & 0 & 0 \\ 
0 & 1 & 0 & 0 & 0 & 0 & 0 & 0 & 0 & 0 \\ 
1 & 0 & 0 & 0 & 0 & 0 & 0 & 0 & 0 & 0%
\end{array}%
\right] \left[ 
\begin{array}{rrrrr}
-\frac{501}{944} & \frac{53}{472} & \frac{23}{59} & -\frac{189}{944} & \frac{%
759}{944} \\ 
-\frac{1315}{944} & -\frac{27}{472} & \frac{44}{59} & -\frac{123}{944} & 
\frac{1393}{944} \\ 
\frac{907}{1888} & \frac{323}{944} & -\frac{25}{59} & \frac{3}{1888} & -%
\frac{1001}{1888} \\ 
-\frac{253}{236} & -\frac{27}{118} & \frac{58}{59} & -\frac{5}{236} & \frac{%
331}{236} \\ 
\frac{205}{118} & \frac{277}{944} & -\frac{5}{59} & \frac{213}{1888} & -%
\frac{1215}{1888} \\ 
-\frac{191}{472} & \frac{9}{236} & \frac{10}{59} & \frac{41}{472} & \frac{165%
}{472} \\ 
\frac{55}{472} & \frac{11}{236} & -\frac{14}{59} & -\frac{81}{472} & \frac{%
123}{472} \\ 
-\frac{619}{944} & -\frac{171}{472} & \frac{23}{59} & -\frac{307}{944} & 
\frac{1113}{944} \\ 
\frac{53}{118} & -\frac{13}{59} & -\frac{25}{59} & -\frac{33}{118} & \frac{37%
}{118} \\ 
\frac{2699}{1888} & \frac{115}{944} & -\frac{50}{59} & \frac{419}{1888} & -%
\frac{3241}{1888}%
\end{array}%
\right. $

$\ \ \ \ \ \ \ \ \ \ \ \ \ \ \ \ \ \ \ \ \ \ \ \ \ \ \ \ \ \ \ \ \ \ \ \ \ \
\ \ \ \ \ \ \ \ \ \ \ \ \ \ \ \ \ \ \ \ \ \ \ \ \ \ \ \ \ \ \ \ \ \ \ \ \ \
\ \ \ \ \ \ \ \ \ \ \ \ \ \ \left. 
\begin{array}{rrrrr}
\frac{709}{944} & \frac{707}{944} & \frac{561}{944} & \frac{367}{472} & -%
\frac{199}{472} \\ 
\frac{1795}{944} & \frac{1509}{944} & \frac{1399}{944} & \frac{1033}{472} & -%
\frac{609}{472} \\ 
-\frac{251}{1888} & -\frac{221}{1888} & -\frac{863}{1888} & -\frac{313}{944}
& \frac{153}{944} \\ 
\frac{261}{236} & \frac{211}{236} & \frac{337}{236} & \frac{207}{118} & -%
\frac{237}{118} \\ 
-\frac{829}{1888} & -\frac{587}{1888} & -\frac{857}{1888} & -\frac{511}{944}
& \frac{479}{944} \\ 
\frac{31}{472} & -\frac{31}{472} & \frac{163}{472} & \frac{49}{236} & -\frac{%
33}{236} \\ 
\frac{169}{472} & \frac{303}{472} & \frac{173}{472} & \frac{191}{236} & -%
\frac{119}{236} \\ 
\frac{827}{944} & \frac{589}{944} & \frac{1151}{944} & \frac{721}{472} & -%
\frac{553}{472} \\ 
\frac{47}{118} & \frac{71}{118} & \frac{53}{118} & \frac{21}{59} & -\frac{31%
}{59} \\ 
-\frac{2331}{1888} & -\frac{1917}{1888} & -\frac{2847}{1888} & -\frac{1865}{%
944} & \frac{1545}{944}%
\end{array}%
\right] $

$\ \ \ \ \ \ \ \ \ \ \ \ \ \ \ =\left[ 
\begin{array}{cccccccccc}
\frac{2699}{1888} & \frac{115}{944} & -\frac{50}{59} & \frac{419}{1888} & -%
\frac{3241}{1888} & -\frac{2331}{1888} & -\frac{1917}{1888} & -\frac{2847}{%
1888} & -\frac{1865}{944} & \frac{1545}{944} \\ 
\frac{53}{118} & -\frac{13}{59} & -\frac{25}{59} & -\frac{33}{118} & \frac{37%
}{118} & \frac{47}{118} & \frac{71}{118} & \frac{53}{118} & \frac{21}{59} & -%
\frac{31}{59} \\ 
-\frac{619}{944} & -\frac{171}{472} & \frac{23}{59} & -\frac{307}{944} & 
\frac{1113}{944} & \frac{827}{944} & \frac{589}{944} & \frac{1151}{944} & 
\frac{721}{472} & -\frac{553}{472} \\ 
\frac{55}{472} & \frac{11}{236} & -\frac{14}{59} & -\frac{81}{472} & \frac{%
123}{472} & \frac{169}{472} & \frac{303}{472} & \frac{173}{472} & \frac{191}{%
236} & -\frac{119}{236} \\ 
-\frac{191}{472} & \frac{9}{236} & \frac{10}{59} & \frac{41}{472} & \frac{165%
}{472} & \frac{31}{472} & -\frac{31}{472} & \frac{163}{472} & \frac{49}{236}
& -\frac{33}{236} \\ 
\frac{205}{118} & \frac{277}{944} & -\frac{5}{59} & \frac{213}{1888} & -%
\frac{1215}{1888} & -\frac{829}{1888} & -\frac{587}{1888} & -\frac{857}{1888}
& -\frac{511}{944} & \frac{479}{944} \\ 
-\frac{253}{236} & -\frac{27}{118} & \frac{58}{59} & -\frac{5}{236} & \frac{%
331}{236} & \frac{261}{236} & \frac{211}{236} & \frac{337}{236} & \frac{207}{%
118} & -\frac{237}{118} \\ 
\frac{907}{1888} & \frac{323}{944} & -\frac{25}{59} & \frac{3}{1888} & -%
\frac{1001}{1888} & -\frac{251}{1888} & -\frac{221}{1888} & -\frac{863}{1888}
& -\frac{313}{944} & \frac{153}{944} \\ 
-\frac{1315}{944} & -\frac{27}{472} & \frac{44}{59} & -\frac{123}{944} & 
\frac{1393}{944} & \frac{1795}{944} & \frac{1509}{944} & \frac{1399}{944} & 
\frac{1033}{472} & -\frac{609}{472} \\ 
-\frac{501}{944} & \frac{53}{472} & \frac{23}{59} & -\frac{189}{944} & \frac{%
759}{944} & \frac{709}{944} & \frac{707}{944} & \frac{561}{944} & \frac{367}{%
472} & -\frac{199}{472}%
\end{array}%
\right] $
\end{itemize}

\begin{example}
Solve the periodic pentadiagonal system given as follows:%
\begin{equation*}
\underset{M}{\underbrace{\left[ 
\begin{array}{cccccc}
1 & 2 & -1 & 0 & 0 & 1 \\ 
2 & -1 & -3 & 1 & 0 & 0 \\ 
1 & 1 & -1 & 1 & 2 & 0 \\ 
0 & 2 & 1 & 1 & -1 & -2 \\ 
0 & 0 & -1 & -2 & 1 & 3 \\ 
1 & 0 & 0 & 1 & 1 & 1%
\end{array}%
\right] }}\left[ 
\begin{array}{c}
x_{1} \\ 
x_{2} \\ 
x_{3} \\ 
x_{4} \\ 
x_{5} \\ 
x_{6}%
\end{array}%
\right] =\left[ 
\begin{array}{c}
3 \\ 
-1 \\ 
4 \\ 
1 \\ 
1 \\ 
4%
\end{array}%
\right]
\end{equation*}%
We apply the Algorithm2:
\end{example}

\begin{itemize}
\item For $i=1,2,3,4$ $\ \ \ \ u_{i,i+1}=(u_{12},u_{23},u_{34},u_{45})=(2,-1,%
\frac{4}{5},-7)$

\item For $i=1,2,3$ $\ \ \ \ u_{i,i+2}=(u_{13},u_{24},u_{35})=(-1,1,2)$

\item For $i=2,3,4,5\ \ \ \ \ l_{i,i-1}=(l_{21},l_{32},l_{43},l_{54})=(2,%
\frac{1}{5},3,-2)$

\item For $i=3,4,5\ \ \ \ \ l_{i,i-2}=(l_{31},l_{42},l_{53})=(1,-\frac{2}{5}%
,-5)$

\item For $i=1,2,3,4,5,6\ \ \ \ \
u_{ii}=(u_{11},u_{22},u_{33},u_{44},u_{55},u_{66})=(1,-5,\frac{1}{5},-1,-3,-%
\frac{14}{3})$

\item For $i=1,2,3,4,5\ \ \ \ \
u_{i,6}=(u_{16},u_{26},u_{36},u_{46},u_{56})=(1,-2,-\frac{3}{5},-1,-2)$

\item For$\ i=1,2,3,4,5\ \ \ \ \
l_{6,i}=(l_{61},l_{62},l_{63},l_{64},l_{65})=(1,\frac{2}{5},7,5,-\frac{22}{3}%
)$

\item $\det (M)=14$

\item $(z_{1},z_{2},z_{3},z_{4},z_{5},z_{6})=(3,-7,\frac{12}{5},-9,-5,-\frac{%
4}{3})$

\item $(x_{1},x_{2},x_{3},x_{4},x_{5},x_{6})=(1,1,1,1,1,1)$
\end{itemize}

\end{document}